\documentclass[a4paper]{article}
\usepackage{ifthen}
\def\style{preprint}
\ifthenelse{\equal{\style}{preprint}}{
\usepackage{dmo_preprint}
}{}
\usepackage{amsthm,amsmath,amsfonts,graphicx}
\usepackage[T1]{fontenc}
\usepackage{amssymb}

\DeclareMathOperator{\rk}{rk}

\DeclareMathOperator{\sgn}{sgn}

\newtheorem{theorem}{Theorem}[section]
\newtheorem*{theorem*}{Theorem}
\newtheorem{lemma}[theorem]{Lemma}
\newtheorem{corollary}[theorem]{Corollary}

\newtheorem{definition}[theorem]{Definition}
\newtheorem{proposition}[theorem]{Proposition}
\newtheorem{observation}[theorem]{Observation}
\newtheorem{example}[theorem]{Example}

\title{A Trivariate Dichromate Polynomial for Digraphs}
\author{Winfried Hochst\"attler, Johanna Wiehe \\ FernUniversit\"at in Hagen, Germany}
\date{}
\begin{document}

\ifthenelse{\equal{\style}{preprint}}{ \DMOmathsubject{05C20,05C31,05C21,05C15}
  \DMOkeywords{flow polynomial, chromatic polynomial, dichromatic number, face lattice, oriented matroids} \DMOtitle{075.23}{A Trivariate Dichromate Polynomial for Digraphs} {Winfried Hochst\"attler, Johanna
    Wiehe}{winfried.hochstaettler@fernuni-hagen.de,
    johanna.wiehe@fernuni-hagen.de} }{}
\maketitle

\begin{abstract}
We define a trivariate polynomial combining the NL-coflow and the NL-flow polynomial, which build a dual pair counting acyclic colorings of directed graphs, in the more general setting of regular oriented matroids.
\end{abstract}

\section{Introduction}

In 1954 Tutte introduced a bivariate polynomial of an undirected graph $G$ and called it the \emph{dichromate} of $G$ \cite{tuttecon}. Nowadays better known as the \emph{Tutte polynomial} it features not only a variety of properties and applications, but also specializes to many graph-theoretic polynomials. Two of them, the chromatic and the flow polynomial, counting proper colorings and nowhere-zero flows, build a pair of dual polynomials in the sense that one polynomial becomes the other one by taking the dual graph. \\
Regarding directed graphs, or digraphs for short, acyclic colorings are a natural generalization of proper colorings. A digraph is \emph{acyclically colorable} if no color class contains a directed cycle. This concept is due to Neumann-Lara \cite{neumannlara}.\\
In \cite{hochstarticle} a flow theory for digraphs transferring Tutte's theory of nowhere-zero flows to directed graphs has been developed and amplified in \cite{althowie} and \cite{howie}, where the authors introduce a pair of dual polynomials, counting acyclic colorings of a digraph and the dual equivalent called NL-flows. \\
In order to combine these two polynomials we will leave the setting of digraphs and enter the world of oriented matroids. This more general scenery provides a plethora of useful techniques as well as a common foundation upon which our new polynomial is built. This foundation is due to a construction of Brylawski and Ziegler \cite{dualpairs} representing a dual pair of oriented matroids as complementary minors. \\
Our notation is fairly standard and follows the book of Bj\"orner et al.\ \cite{BLSW} if not explicitely defined. \\

In \cite{howie} we found the following representation of the NL-coflow polynomial counting acyclic colorings in a digraph $D=(V,A)$. A digraph is \emph{totally cyclic}, if every component is strongly connected.\\

\begin{definition}
Let $\mu_Q$ be the M\"obius function of $(Q,\subseteq)$ with 
$$Q:=\lbrace B \subseteq A : D[B] \text{ is a totally cyclic subdigraph of } D \rbrace.$$ 
Then
$$\psi_{NL}^D(x)=\sum_{B \in Q} \mu_Q(\mathbf{0},B)x^{rk(A/B)}$$
is called the \emph{NL-coflow polynomial} of $D$, where, for $Y\subseteq A$, $rk(Y)$ is the rank of the incidence matrix of $D[Y]$, which equals $|V(Y)|-c(Y)$, i.e.\ the number of vertices minus the number of (weakly) connected components of $D[Y]$.
\end{definition}
Recall that in our definition of contraction (see \cite{bondy}) no additional arcs (elements) are removed, i.e. parallel arcs and loops can occur. This holds for both the graphic and the matroid contraction.\\
We will now define this polynomial in the more general setting of (regular) oriented matroids. Note, that all of our results also work in the non-regular case. Since we are not aware of a meaningful interpretation in this case, all our matroids will be regular, if not explicitely pointed out.\\
Let $M$ be an oriented matroid on a (finite) groundset $E$. Let $X,Y\in \lbrace 0,+,-\rbrace^{|E|}$ be two signed sets. The \emph{composition} of $X$ and $Y$, denoted by $X\circ Y$, is the signed set $Z\in \lbrace 0,+,-\rbrace^{|E|}$, where $Z_i=X_i$, if $X_i\neq 0$ and $Z_i=Y_i$, otherwise. The covectors of $M$, i.e.\ compositions of (signed) cocircuits, together with the partial order $0\leq +$ and $0\leq -$ form the \emph{face lattice} $\mathcal{L}$ of $M$ with minimal element $\mathbf{0}:=(0,\ldots,0)$. Since the NL-flow polynomial (see \cite{althowie}) only considers directed cuts, we are only interested in the nonnegative part of $\mathcal{L}$, which we denote by $\mathcal{L}_{+}:=\mathcal{L}\cap \lbrace 0,+\rbrace^E$. By $\mathcal{L}^\ast$ we denote the face lattice of the dual $M^\ast$. Again, we are only interested in the nonnegative part $\mathcal{L}^\ast_+$ which in the graphic case corresponds to the set of totally cyclic subdigraphs partially ordered by inclusion. 
Let $\mu$ and $\mu^\ast$ denote the M\"obius function of $\mathcal{L}_+$ and $\mathcal{L}^\ast_+$, respectively. By $\rk$ and $\rk^\ast$ we denote the rank and corank of the respective matroid (minor) and by $\underline{X}$ we denote the support of the covector $X$. \\
Now we can define the \emph{NL-coflow polynomial} of an oriented matroid $M$ as
$$ \psi_{NL}^M(x):=\sum_{X \in \mathcal{L}^\ast_{+}}\mu^\ast(\mathbf{0},X)x^{\rk(M/\underline{X})}.$$
Dually we define the \emph{NL-flow polynomial} of $M$ as
$$ \phi_{NL}^M(x):= \sum_{X \in \mathcal{L}_{+}}\mu(\mathbf{0},X)x^{\rk^\ast(M\setminus \underline{X})}.$$ 
It is easy to see that both coflow polynomials coincide in the graphic case. Our new definition of the NL-flow polynomial also coincides with the graphic one in \cite{althowie} since $\rk^\ast(Y)=|Y|-\rk(Y)$ holds for any minor, in particular for $Y:=M\setminus B, B \in \mathcal{C}$\footnote{In \cite{althowie} the NL-flow polynomial of a digraph $D=(V,A)$ is defined on the poset $(\mathcal{C},\supseteq)$ with $
\mathcal{C}:= \lbrace A\setminus C : \exists\; C_1,\ldots,C_r\textrm{ directed cuts s.t. } C=\bigcup\limits_{i=1...r} C_i \rbrace.$
}.\\
In particular, the M\"obius function simply alternates with respect to the rank of $\mathcal{L}_+$ as the following shows.\\

\begin{lemma}\label{coneOM}
Let $M$ be an oriented matroid. Then $\mathcal{L}_+$ ordered componentwise by $0\leq +$ is a graded lattice with rank function $\rk_{\mathcal{L}_+}=\rk_{\mathcal{L}}$ and its M\"obius function alternates in the following fashion:
$$ \mu_{\mathcal{L}_+}(\mathbf{0},X)=(-1)^{\rk_{\mathcal{L}_+}(X)},$$
for all $X\in \mathcal{L}_+$.
\end{lemma}

\begin{proof}
By Theorem 4.1.14 in \cite{BLSW} the face lattice $\mathcal{L}$ is graded. Let $X_+$ denote the covector in $\mathcal{L}$ with maximal number of $+$-entries. Note that $X_+$ is uniquely determined. Then $\mathcal{L}_+=[\mathbf{0},X_+]$ forms a lattice itself, inheriting the rank function of $\mathcal{L}$. Corollary 4.3.8 in \cite{BLSW} completes the proof.
\end{proof}

\section{Setting}

Since our polynomials are defined on different face lattices we have to find a common lattice including both. In \cite{dualpairs} Brylawski and Ziegler give the following beautiful construction which provides the desired lattice. \\
Let $M$ be an oriented matroid on the groundset $E=\lbrace 1,\ldots, n\rbrace$ with rank $r$ and $M^\ast$ its dual. Suppose that $\mathcal{B}:=\lbrace 1,\ldots,r\rbrace$ is a basis of $M$ and $\lbrace r+1,\ldots,n\rbrace$ is the corresponding basis of $M^\ast$. Furthermore set $E_1:=\mathcal{B}, E_2:=E\setminus \mathcal{B}$ and
$$\hat{E}:=E_1\cup E_2 \cup A \cup B = E\cup A \cup B$$ 
with $A:=\lbrace n+1,\ldots ,n+r \rbrace$ and $B:=\lbrace n+r+1,\ldots, 2n \rbrace$ and let $M_1$ be the oriented matroid on $\hat{E}$, that is obtained by extending $M$ by elements $n+i$ that are parallel to the elements $i$ for $1\leq i \leq r$ and that are loops for $r+1\leq i \leq n$. Similarly, let $M_2$ be the oriented matroid on $\hat{E}$ that is obtained by extending $M^\ast$ by elements $n+i$ that are loops for $1\leq i \leq r$ and that are parallel to the elements $i$ for $r+1\leq i \leq n$. Then $M_1$ has rank $r$ and $M_2$ has rank $n-r$. Their union (see chapter 7.6 in \cite{BLSW}) 
$$\hat{M}:=M_1 \cup M_2$$ 
has rank $n$. Note that the construction of the oriented matroid union highly depends on the choice of the basis $\mathcal{B}$. Due to Theorem 2 in \cite{dualpairs} we have 
$$ \hat{M}\setminus A/B=M \textrm{ and } \hat{M}/A\setminus B= M^\ast.$$
In the case where $M$ is realizable, $\hat{M}$ is also realizable. Namely, if $M$ can be represented by $\begin{pmatrix}
I_r & C
\end{pmatrix}$, where $I_r$ denotes the identity matrix of rank $r$, then $M_1$ and $M_2$ are represented by $\begin{pmatrix}
I_r&C&I_r&0
\end{pmatrix}$ and $\begin{pmatrix}
-C^\top & I_{n-r}& 0& I_{n-r}
\end{pmatrix}$, respectively. Now let $\begin{pmatrix}
-C^\top & I_{n-r}& 0 &I_{n-r}
\end{pmatrix}^\epsilon$ be the matrix obtained by multiplying the $i-$th column by $\epsilon^{2n-i}$ for all \mbox{$i\in  \lbrace 1,\ldots ,2n\rbrace$} and $\epsilon > 0$ sufficiently small. Then the combined matrix 
$$\begin{pmatrix}
I_r  & C& I_r& 0\\
(-C^\top &I_{n-r} &0 &I_{n-r} )^\epsilon
\end{pmatrix}$$
represents $\hat{M}$ (see Proposition 8.2.7 of \cite{BLSW} and \cite{dualpairs}). Note that even if $M$ and $M^\ast$ are regular, this might not be true for $\hat{M}$. However, the face lattice of $\hat{M}$, which we will denote by $\hat{\mathcal{L}}$, will serve our purpose.\\ 
The upcoming analysis yields that this construction is not only elegant, it also preserves the natural separation of the totally cyclic part and the acyclic part (see Corollary 3.4.6 in \cite{BLSW}). In the following we will find a characterization of the covectors of $M$ and its dual in this supermatroid $\hat{M}$. But let us first define the new polynomial in order to motivate the upcoming rather technical analysis of $\hat{M}$.

\section{The Dichromate}

\subsection{Definition}

In this setting we are able to define a new polynomial in three variables which generalizes both, the NL-flow and the NL-coflow polynomial, meaning that both can be found as evaluations of the new one. In order to switch between the NL-flow and the NL-coflow polynomial we use two of the three variables as some kind of toggle. We will prove that whenever the support of a covector of $\hat{M}$ is non-empty on $A$ (or on $B$ resp.), this covector cannot correspond to one of $M^\ast$ (or $M$ resp.) and will be rejected if $y$ or $z$ are evaluated at zero. In Lemma \ref{restrict} we prove that covectors that correspond neither to a covector of $M$ nor to one of $M^\ast$ will also be rejected, since they have non-empty support on $A$ as well as on $B$. This is why we can define our polynomial on the whole face lattice $\hat{\mathcal{L}}_{+}$.\\

\begin{definition}\label{Def.dichromate}
Let $M$ be a regular, oriented matroid on a finite groundset $E$, $\mathcal{B}$ the basis of $M$ chosen to construct $\hat{M}$ and $\hat{\mu}$ the M\"obius function of the face lattice of $\hat{M}$. Then we define
$$\Omega_{NL}^{M,\mathcal{B}}(x,y,z):= \sum_{X \in \hat{\mathcal{L}}_{+}} \hat{\mu}(\mathbf{0},X) x^{\rk_{\hat{\mathcal{L}}}(X)+|E\setminus (\underline{X}\cap E)|}y^{|\underline{X}\cap A|}z^{|\underline{X}\cap B|},$$
which we call the \emph{dichromate} of a digraph representing $M$ in the graphic case.
\end{definition}

Note that,
$$ \rk^\ast(\hat{M}\setminus \underline{X})=\rk_{\hat{\mathcal{L}}}(X)+|\hat{E}\setminus \underline{X}|-\rk(\hat{M}).\footnote{see (\ref{rank-cor.}) and Lemma \ref{latticerank}(ii)}$$
Therefore, the dichromate of $M$ is closely related to the NL-flow polynomial of $\hat{M}$ in the $x$-variables. The $y$-variables control the elements parallel to the chosen basis while the $z$-variables control the cobasis.

\subsection{Properties}

\subsubsection{Evaluations}

A first, simple observation is the following.\\

\begin{observation}\label{111}
Let $M$ be a regular, oriented matroid of rank $r$. Then
$$\Omega_{NL}^{M,\mathcal{B}}(1,1,1)=0$$
for any basis $\mathcal{B}$ of $M$.
\end{observation}
\begin{proof}
In Lemma \ref{coneOM} we have seen that the M\"obius function of $\hat{\mathcal{L}}_+$ alternates with respect to the rank. Using the \emph{Euler-Poincar\'{e} formula} (see e.g.\ Corollary 4.6.11 in \cite{BLSW}) we find:
$$\sum_{X\in \hat{\mathcal{L}}_+}\hat{\mu}(\mathbf{0},X)= \sum_{X\in \hat{\mathcal{L}}_+}(-1)^{\rk_{\hat{\mathcal{L}}_+}(X)}=(-1)^{\rk_{\hat{\mathcal{L}}_+}(\mathbf{0})}-(1+(-1)^{r-1})+(-1)^{r}=0.$$
\end{proof}

Our main goal is to prove the following theorems.\\

\begin{theorem}\label{trivariatecoflow}
Let $M$ be a regular, oriented matroid on $E$ with $|E|=n$ and let $r$ be its rank. Then
$$\Omega_{NL}^{M,\mathcal{B}}(x,0,1)=x^{n-r}\cdot \psi_{NL}^M(x)$$
for any basis $\mathcal{B}$ of $M$.\\
\end{theorem}

\begin{theorem}\label{trivariateflow}
Let $M$ be a regular, oriented matroid on $E$ with $|E|=n$ and let $r$ be its rank. Then
$$ \Omega_{NL}^{M,\mathcal{B}}(x,1,0)=x^r\cdot\phi_{NL}^M(x)$$
for any basis $\mathcal{B}$ of $M$.
\end{theorem}

\paragraph{Cocircuits and Covectors}\label{sec:cocircuits}

Given a cocircuit $D$ in $M$ or in $M^\ast$ we find a corresponding cocircuit $\hat{D}$ in $\hat{M}$ such that $\underline{D}\subseteq \underline{\hat{D}}$. Furthermore we will find that, given $D^-=\emptyset$, then also $\hat{D}^-=\emptyset$ holds. Due to the construction of $\hat{M}$ we will first extend $D$ to a cocircuit in $M_1$, which then is already a cocircuit in $\hat{M}$. For the proof we will first look at the underlying unoriented matroid and then compute the signatures in a second step. For an element $a \in \hat{E}$ we write $a||b$, if there exists an $i\in [n]$ and some $b\in A\cup B$, such that $a=i$ and $b=n+i$.\\

\begin{lemma}\label{undirected1}
Let $D$ be a cocircuit in $\underline{M}$ and set $D_1:= \lbrace a \in A : e||a, e \in D\cap{E_1}\rbrace.$ Then $\hat{D}:= D\cup D_1$ is a cocircuit in $\underline{M_1}$. 
\end{lemma}

\begin{proof}
First we prove that $\hat{D}$ meets every basis of $\underline{M_1}$.\\
Let $b$ be a basis of $\underline{M_1}$. Since all elements in $B$ are loops $b \subseteq E\cup A$ has to hold. If $b \subseteq E$ then $b$ is a basis of $\underline{M}$ and $\hat{D}\cap b = (D\cap b)\cup (D_1 \cap b) \neq \emptyset$. If $b \cap A\neq \emptyset$ we find a basis
$$ b^\prime= (b\cap E) \cup \lbrace e \in E_1 : e||f, f \in b\cap A\rbrace\subseteq E$$
in $\underline{M}$. Hence $D \cap b^\prime\neq \emptyset$ and in particular $\hat{D}\cap b^\prime \neq \emptyset$. So either $\hat{D}\cap (b\cap E)\neq \emptyset$ or $b=A$ but then $D_1\cap b\neq \emptyset$ and in any case $\hat{D}\cap b\neq \emptyset$. Note, that $D_1\neq \emptyset$, otherwise $D\cap E_1=\emptyset$ which cannot be since $E_1$ is a basis of $\underline{M}$. \\
For the minimality of $\hat{D}$, let $d \in \hat{D}$. If $d \in D$ we find a basis $b$ in $\underline{M}$ with $(D\setminus d )\cap b =\emptyset$. Then also $(\hat{D}\setminus d)\cap b=\emptyset$ since $D_1\cap b=\emptyset$. Otherwise, if $d \in D_1$, there exists an $f \in D\cap E_1$ with $f||d$. Thus there exists a basis $b$ in $\underline{M}$, which is also a basis of $\underline{M_1}$, with $(D\setminus f)\cap b =\emptyset$. Due to the definition of $D_1$ we also have that $(D_1 \setminus d)\cap b =\emptyset$. Then, by basis exchange
$$b_1:=(b\setminus f) \cup d$$
is a basis of $\underline{M_1}$ with $(\hat{D}\setminus d)\cap b_1=\emptyset$.\\
\end{proof}

\begin{lemma}\label{directed1}
If $D=(D^+,D^-)$ is a signed cocircuit in $M$ with $D^-=\emptyset$, then $\hat{D}:=(D^+\cup D_1 ,\emptyset)$ is a signed cocircuit in $M_1$. 
\end{lemma}

\begin{proof}
With $\chi_1$ and $\chi_M$ denote the chirotope of $M_1$ and $M$ respectively. \\
Let $e,f \in \underline{D}, e\neq f$ and $(x_2,\ldots,x_r)$ be any ordered basis of the hyperplane $E\setminus \underline{D}$ in $M$. Furthermore set
$$ \sigma_M(e,f):= \chi_M(e,x_2,\ldots,x_r)\chi_M(f,x_2,\ldots,x_r) \in \lbrace 1,-1\rbrace.$$
Then, due to Lemma 3.5.8 in \cite{BLSW} and since $D^-=\lbrace f \in \underline{D}\setminus e : \sigma_M(e,f)=-1\rbrace= \emptyset$ by assumption, we have $\sigma_M(e,f)=1$ for all $f\in \underline{D}\setminus e$.\\
Now, let $e,f \in \underline{\hat{D}}\cap \underline{D}, e\neq f$ and let $X=(x_2,\ldots ,x_r)$ be an ordered basis of the hyperplane $\hat{E}\setminus \underline{\hat{D}}$ in $M_1$. Let $k$ be the first index such that $x_i \in A$ for all $i\geq k$ and $x_i \in E$ for all $i <k$. Since all elements in $A$ are parallel to the elements in $E_1$ in $M_1$ we find a basis $X^\prime=(x_2^\prime,\ldots , x_r^\prime)=(x_2,\ldots,x_{k-1},y_{k},\ldots ,y_r)$ of the hyperplane $E\setminus \underline{D}$ in $M$ by mapping all elements $x_i$ from $X\cap A$ to their corresponding parallels $y_i$ in $E_1$ which cannot be in $X$ since this is a basis. Now let $\tau$ be a permutation of the elements in $X^\prime$, such that $\tau(X^\prime)=(x_{\tau(2)}^\prime,\ldots, x_{\tau(r)}^\prime)$ is ordered in $E$. Then we find
\begin{align*}
\sigma_1(e,f)&= \chi_1(e,x_2,\ldots, x_r)\cdot \chi_1(f,x_2,\ldots ,x_r)\\
&=\chi_1(e,x_2,\ldots ,x_{k-1},y_k,\ldots ,y_r)\cdot \chi_1(f,x_2,\ldots x_{k-1},y_k,\ldots y_r)\\
&=\chi_1(e,x_{\tau(2)}^\prime,\ldots,x_{\tau(r)}^\prime)\cdot\chi_1(f,x_{\tau(2)}^\prime,\ldots, x_{\tau(r)}^\prime)\cdot \sgn(\tau)^2
=\sigma_M(e,f)=1.
\end{align*}
Thus $f \in \hat{D}^+$ for all $f \in \underline{\hat{D}}\cap \underline{D}$.
Otherwise, if $f \in \underline{\hat{D}}\cap \underline{D_1}$, then there exists $g \in E_1, g||f$ with $g \in \underline{D}$. Similarly as above we find $\sigma_1(e,g)=\sigma_M(e,g)=1$
and so also $f \in \hat{D}^+$ has to hold for all $f \in \underline{\hat{D}}\cap \underline{D_1}$. As a result $\hat{D}^-=\emptyset$.
\end{proof}
We are left to prove that $\hat{D}$ is also a (signed) cocircuit in $\hat{M}$. Again, we will first take a look at the underlying unoriented case, where the oriented matroid union becomes the usual matroid union. Let $\mathcal{I}_1$ and $\mathcal{I}_2$ be the independent sets in $\underline{M_1}$ and $\underline{M_2}$ respectively. Then $\underline{\hat{M}}=(\hat{E},\mathcal{I})$, where 
$$\mathcal{I}=\lbrace I_1 \cup I_2 : I_1 \in \mathcal{I}_1 \text{ and } I_2\in \mathcal{I}_2\rbrace$$
are the independent sets of $\underline{\hat{M}}$. As an immediate result every basis $b$ of $\underline{\hat{M}}$ can be written as $b=b_1\cup b_2$, where $b_1$ is a basis of $\underline{M_1}$ and $b_2$ is a basis of $\underline{M_2}$. Note that since $|b_1|=r$ and $|b_2|=n-r$ the union $b_1\cup b_2$ must be disjoint.\\

\begin{lemma}\label{undirected2}
Let $D$ be a cocircuit in $\underline{M_1}$. Then $D$ is also a cocircuit in $\underline{\hat{M}}$.
\end{lemma}

\begin{proof}
Let $b=b_1\cup b_2$ be a basis of $\underline{\hat{M}}$, where $b_1$ is a basis of $\underline{M_1}$ and $b_2$ is a basis of $\underline{M_2}$. Since $D$ is a cocircuit in $\underline{M_1}$, in particular $D\cap b_1\neq \emptyset$ and so $D\cap b\neq \emptyset$ has to hold.\\
For the minimality let $d \in D$. Since $D$ is a cocircuit in $\underline{M_1}$ we find $(D\setminus d)\cap b_1=\emptyset$ for some basis $b_1$ in $\underline{M_1}$. Then $b:=b_1\cup B$ is a basis of $\underline{\hat{M}}$ since $B$ is a basis of $\underline{M_2}$ and $B\cap b_1=\emptyset$ since all elements in $B$ are loops in $\underline{M_1}$. As a result $ (D\setminus d)\cap b=\emptyset.$\\
\end{proof}

\begin{lemma}\label{directed2}
If $D=(D^+,D^-)$ is a signed cocircuit in $M_1$, then $D$ is a signed cocircuit in $\hat{M}$.
\end{lemma}

\begin{proof}
Let $e,f\in \underline{D}$, $e\neq f$ and $X:=(x_2,\ldots,x_n)$ be an ordered basis of the hyperplane $\hat{E}\setminus \underline{D}$ in $\hat{M}$. Since $\underline{D}\cap B =\emptyset$ we can choose $X$ such that $B\subseteq X$. Hence $(x_2,\ldots,x_r)$ is an ordered basis of the hyperplane $\hat{E}\setminus \underline{D}$ in $M_1$. Now let $\sigma$ and $\tau$ be permutations of the elements in $x:=(e,x_2,\ldots,x_r)$ and in $y:=(f,x_2,\ldots,x_r)$ respectively, such that $(\sigma(x),B)$ and $(\sigma(y),B)$ are both lexicographically minimal ordered. This is always possible since no element in $B$ can be part of a basis of $M_1$. Then by Proposition 7.6.4 in \cite{BLSW}
\begin{flalign*}
\hat{\sigma}(e,f)&=\hat{\chi}(e,x_2,\ldots,x_r,x_{r+1},\ldots,x_n)\cdot \hat{\chi}(f,x_2,\ldots,x_r,x_{r+1},\ldots,x_n)\\
&=\hat{\chi}(\sigma(x),B)\cdot \sgn(\sigma)\cdot \hat{\chi}(\tau(y),B)\cdot\sgn(\tau)\\
&=\chi_1(\sigma(x))\cdot \chi_2(B)\cdot \sgn(\sigma)\cdot \chi_1(\tau(y))\cdot \chi_2(B)\cdot \sgn(\tau)\\
&=\chi_1(e,x_2,\ldots,x_r)\cdot \chi_1(f,x_2,\ldots,x_r)\cdot \sgn^2(\sigma)\cdot \sgn^2(\tau)\\
&=\sigma_1(e,f).
\end{flalign*}
Lemma 3.5.8 in \cite{BLSW} completes the proof.
\end{proof}

Analogously, one can define $D_2:=\lbrace b \in B : e||b, e \in D\cap{E_2} \rbrace$ and prove that if $D=(D^+,\emptyset)$ is a signed cocircuit in $M^\ast$, then $\hat{D}:=(D^+\cup D_2,\emptyset)$ 
is a signed cocircuit in $\hat{M}$. Since covectors are compositions of cocircuits, the results above readily yield:\\

\begin{proposition}
\hfill
\begin{enumerate}
\item[(i)] Let $X$ be a covector in $M$ and $\tilde{A}:=\lbrace a \in A: e||a, e \in \underline{X}\cap{E_1}\rbrace$. Then $\hat{X}:=(X^+\cup\tilde{A},\emptyset)$ is a covector in $\hat{M}$.
\item[(ii)] Let $X$ be a covector in $M^\ast$ and $\tilde{B}:=\lbrace b \in B: e||b, e \in \underline{X}\cap{E_2}\rbrace.$ Then $\hat{X}:=(X^+\cup \tilde{B},\emptyset)$ is a covector in $\hat{M}$. 
\end{enumerate}
\end{proposition}

\paragraph{The Face Lattice of $\hat{M}$}

We have already seen that both the covectors of $M$ and of $M^\ast$ can be found in the face lattice of $\hat{M}$. In the following we will show the converse: Having a covector of $\hat{M}$ of that specific shape we determined in the previous section, its restriction to $E$ corresponds to a covector of $M$ or of $M^\ast$, respectively. Furthermore we will see that the corresponding M\"obius functions coincide. The following lemma will be crucial for both. Here we write $(\underline{\hat{X}} \cap E_1)||(\underline{\hat{X}} \cap A)$, if $x\in E_1$ is in $\underline{\hat{X}}$ if and only if $(x+n)\in A$ is also in $\underline{\hat{X}}$.\\
Due to a recurrent analogy we will only give the proofs of each of the first alternatives.\\ 

\begin{lemma}\label{crucial}
Let $\hat{X}=(\hat{X}^+,\emptyset)$ be a signed cocircuit of $\hat{M}$ with 
$\underline{\hat{X}}\cap B =\emptyset$ \mbox{($\underline{\hat{X}} \cap A =\emptyset$).} Then $(\underline{\hat{X}} \cap {E_1})||(\underline{\hat{X}} \cap A)$ (resp. $(\underline{\hat{X}} \cap {E_2}) ||(\underline{\hat{X}} \cap B)$).
\end{lemma}

\begin{proof}
Let $x\in \underline{\hat{X}}\cap {E_1}$. Since $\hat{X}$ is a cocircuit of $\hat{M}$ there exists a basis $b=b_1\cup b_2$ of $\hat{M}$ such that $(\underline{\hat{X}}\setminus x)\cap b=\emptyset$, where $b_1$ is a basis of $M_1$ and $b_2$ of $M_2$. By assumption $\underline{\hat{X}}\cap B=\emptyset$, thus $\underline{\hat{X}}\cap b_1$ must be nonempty. It follows that $(\underline{\hat{X}}\setminus x)\cap b_1=\emptyset$, $x\in b_1$ and $\underline{\hat{X}}\cap b_2=\emptyset$. Let $y\in A$ such that $x||y$. By basis exchange we obtain a new basis 
$$ b_1^\prime := (b_1\setminus x) \cup y$$
of $M_1$. Note, that $y$ cannot be in any basis of $M_2$, therefore $b_1^\prime \cup b_2$ is a basis of $\hat{M}$. Since $\hat{X}$ is a cocircuit also 
$$\emptyset \neq \underline{\hat{X}}\cap b^\prime=(\underline{\hat{X}}\cap b_1^\prime)\cup (\underline{\hat{X}}\cap b_2)=(\underline{\hat{X}}\cap (b_1\setminus x))\cup (\underline{\hat{X}}\cap y)=\underline{\hat{X}}\cap y$$ 
holds. Hence $y\in \underline{\hat{X}}\cap A$. The other direction can be proven similarly.
\end{proof}

For a signed cocircuit $\hat{X}=(\hat{X}^+,\hat{X}^-)$ of $\hat{M}$ we set $\hat{X}\cap E:=(\hat{X}^+\cap E,\hat{X}^-\cap E)$.\\

\begin{lemma}\label{restrict}
Let $\hat{X}=(\hat{X}^+,\emptyset)$ be a signed cocircuit of $\hat{M}$.
\begin{itemize}
\item[(i)] If $\underline{\hat{X}}\cap B=\emptyset$, then $X:=\hat{X}\cap E$ is a signed cocircuit of $M$ and $X^-=\emptyset$.
\item[(ii)] If $\underline{\hat{X}}\cap A=\emptyset$, then $X:=\hat{X}\cap E$ is a signed cocircuit of $M^\ast$ and $X^-=\emptyset$.
\end{itemize}
\end{lemma}

\begin{proof}
Let $B_1$ be a basis of $M$. Then $B_1 \cup B$ is a basis of $\hat{M}$. Since $\hat{X}$ is a cocircuit of $\hat{M}$ it meets every basis, in particular we find $ \underline{\hat{X}} \cap (B_1\cup B) \neq \emptyset.$ As we have $\underline{\hat{X}}\cap B=\emptyset$ it follows immediately that $\hat{X}$ meets every basis $B_1$ of $M$. Since $B_1 \subseteq E$ also $\underline{X}\cap B_1 \neq \emptyset$ has to hold.\\
For the minimality of $X$ in $M$ remove $x\in E$ from $\hat{X}$, then there exists a basis $(b_1 \cup b_2)$ of $\hat{M}$ such that $(\underline{\hat{X}}\setminus x) \cap (b_1\cup b_2) =\emptyset.$ As a result also $(\underline{\hat{X}}\setminus x) \cap b_1=\emptyset$ holds for some basis $b_1$ of $M_1$. Now let $b$ be the basis of $M$, where all the elements of $b_1 \cap A$ are replaced by their parallels in $E_1$, i.e.
$$ b:= (b_1\cap E)\cup \lbrace e \in E_1 : e||f, f\in b_1\cap A\rbrace.$$
Due to Lemma \ref{crucial} $(\underline{\hat{X}}\cap {E_1})||(\underline{\hat{X}}\cap A)$ holds, which yields that $ (\underline{X}\setminus x)\cap b=\emptyset.$\\
In order to determine the signatures of $X$ let $e \in \underline{\hat{X}}$ and $(x_2,\ldots,x_n)$ be an ordered basis of the hyperplane $\hat{E}\setminus \underline{\hat{X}}$ such that $(x_2,\ldots,x_r)$ is a basis of the ordered hyperplane $E\setminus \underline{X}$. Since $\hat{X}^-=\emptyset$ Lemma 3.5.8 in \cite{BLSW} yields that
$$\hat{\sigma}(e,f)=\hat{\chi}(e,x_2,\ldots,x_n)\cdot \hat{\chi}(f,x_2,\ldots,x_n)=1$$ 
for all $f \in \underline{\hat{X}}\setminus e$. As we have $M=\hat{M}\setminus A / B$ we find
\begin{align*}
\sigma_M(e,f)&=\chi_M(e,x_2,\ldots ,x_r) \chi_M(f,x_2,\ldots,x_r)\\
&= \hat{\chi}(e,x_2,\ldots ,x_r,b_1,\ldots, b_{n-r})\hat{\chi}(f,x_2,\ldots ,x_r,b_1,\ldots, b_{n-r})=1,
\end{align*}
where $b_1,\ldots,b_{n-r}$ is the ordered basis of $B$ in $\hat{M}$. 
\end{proof}
Again, the previous lemma generalizes naturally to covectors. Let us now take a look at the corresponding rank functions. By $\rk_{\mathcal{L}}, \rk_{\mathcal{L}^\ast}$ and $\rk_{\hat{\mathcal{L}}}$ we denote the rank functions of the respective face lattices of $M, M^\ast$ and $\hat{M}$.\\

\begin{lemma}\label{rank}
Let $X=(X^+,\emptyset)$ be a covector of $M$ (of $M^\ast$) and let $\hat{X}$ be the corresponding covector in $\hat{M}$. Then $\rk_{\mathcal{L}}(X)=\rk_{\hat{\mathcal{L}}}(\hat{X})$ ($\rk_{\mathcal{L}^\ast}(X)=\rk_{\hat{\mathcal{L}}}(\hat{X})$).
\end{lemma}

\begin{proof}
Since all the covectors $Y$ with $\underline{Y}\subseteq \underline{X}$ have corresponding covectors $\hat{Y}$ with $\underline{\hat{Y}}\subseteq \underline{\hat{X}}$ it is clear that $\rk_{\mathcal{L}}(X)\geq \rk_{\hat{\mathcal{L}}}(\hat{X})$ holds. \\
For the other direction let $\hat{X}=Y\circ \hat{Z}$, where $Y$ is a cocircuit of $\hat{M}$. Then, since $\underline{\hat{X}}\cap B =\emptyset$ also $\underline{\hat{Z}}\cap B=\emptyset$ holds and so, due to Lemma \ref{restrict}, $\underline{\hat{Z}}\cap E$ is a covector of $M$. Inductively we get $\rk_{\mathcal{L}}(X)\leq \rk_{\hat{\mathcal{L}}}(\hat{X})$ completing the proof.
\end{proof}
As an immediate result, also the corresponding M\"obius functions coincide.  
Aside from this we will find a common expression of the exponents of the NL-flow and the NL-coflow polynomial in terms of the rank in the face lattice of $\hat{M}$. In order to do so we will use Corollary 4.1.15 (i) in \cite{BLSW}:
\begin{eqnarray}\label{rank-cor.}
\rk_{\mathcal{L}}(X)=\rk(M)-\rk(M\setminus \underline{X}) \quad \forall X \in \mathcal{L}.
\end{eqnarray}
\begin{lemma}\label{latticerank}
Let $\hat{X} \in \hat{\mathcal{L}}_+$ and $X:=\hat{X}\cap E$.
\begin{itemize}
\item[(i)] If $\underline{\hat{X}}\cap A=\emptyset$, then 
$ \rk(M/\underline{X})=\rk_{\hat{\mathcal{L}}}(\hat{X})+|E\setminus \underline{X}|-(n-r).$
\item[(ii)] If $\underline{\hat{X}}\cap B=\emptyset$, then $ \rk^\ast(M\setminus \underline{X})= \rk_{\hat{\mathcal{L}}}(\hat{X})+|E\setminus \underline{X}|-r.$
\end{itemize}
\end{lemma}

\begin{proof} 
\begin{itemize}
\item[(i)] Due to Lemma \ref{restrict}, $X \in \mathcal{L}^\ast_+$. Dualizing, (\ref{rank-cor.}) and Lemma \ref{rank} yield
\begin{flalign*}
\rk(M/\underline{X})&= \rk^\ast(M^\ast \setminus \underline{X})=|E\setminus \underline{X}|-\rk(M^\ast\setminus \underline{X})\\
&=|E\setminus \underline{X}|+ \rk_{\mathcal{L}^\ast}(X)-\rk(M^\ast)=|E\setminus \underline{X}| + \rk_{\hat{\mathcal{L}}}(\hat{X})-(n-r).
\end{flalign*}
\item[(ii)] Due to Lemma \ref{restrict} $X \in \mathcal{L}_+$. (\ref{rank-cor.}) and Lemma \ref{rank} yield
\begin{flalign*}
\rk^\ast(M\setminus \underline{X})&= |E\setminus \underline{X}|-\rk(M\setminus \underline{X})=|E\setminus \underline{X}|+ \rk_{\mathcal{L}}(X)-\rk(M)\\
&=|E\setminus \underline{X}| + \rk_{\hat{\mathcal{L}}}(\hat{X})-r.
\end{flalign*}
\end{itemize}
\end{proof}

Finally, we can prove Theorem \ref{trivariatecoflow} and Theorem \ref{trivariateflow}.


\begin{proof}[Proof of Theorem \ref{trivariatecoflow}]
Due to the definition it immediately follows that
$$\Omega_{NL}^{M,\mathcal{B}}(x,0,1)=\sum_{\substack{X\in \hat{\mathcal{L}}_{+}\\\underline{X}\cap A =\emptyset}}\hat{\mu}(\mathbf{0},X)x^{\rk_{\hat{\mathcal{L}}}(X)+|E\setminus (\underline{X}\cap E)|}$$
for any basis $\mathcal{B}$ of $M$. Lemma \ref{restrict} (ii) yields that the sum only considers \mbox{$X\cap E\in \mathcal{L}_{+}^\ast$,} since it is a covector of $M^\ast$ with positive entries only. The respective M\"obius functions coincide due to Lemma \ref{rank}. Lemma \ref{latticerank} (i) completes the proof.
\end{proof}
Using Lemmas \ref{restrict} (i), \ref{rank} and \ref{latticerank} (ii), Theorem \ref{trivariateflow} can be proven completely analogously.


\subsubsection{Duality}

Besides those two main properties the trivariate polynomial naturally fulfills a third one, similar to the duality property of the Tutte polynomial $T_M(x,y)$ of a matroid $M$, i.e. (see e.g. \cite{welshTutte})
$$T_M(x,y)=T_{M^\ast}(y,x).$$

\begin{theorem}\label{triduality}
Let $M$ be a regular, oriented matroid on a groundset $E$. Then
$$ \Omega_{NL}^{M,\mathcal{B}}(x,y,z)=\Omega_{NL}^{M^\ast,E\setminus\mathcal{B}}(x,z,y)$$
for any basis $\mathcal{B}$ of $M$.
\end{theorem}

\begin{proof}
Let $E=[n]$ and $r$ denote the rank of $M$. Let $\mathcal{B}=[r]$ and $M$ be represented by $\begin{pmatrix}
I_r& C
\end{pmatrix}$. Then $E\setminus \mathcal{B}=\lbrace r+1,\ldots,n\rbrace$ and $M^\ast$ is represented by $M^\ast:=\begin{pmatrix}
-C^\top& I_{n-r}
\end{pmatrix}$. In order to construct $\hat{(M^\ast)}$ with respect to the basis $E\setminus \mathcal{B}$ we need to convert $M^\ast$ into $\begin{pmatrix}
I_{n-r}-C^\top
\end{pmatrix}$ by applying an isomorphism $\varphi:\ E\longrightarrow E$, 
$$ \varphi(i)=\begin{cases} i+n-r, &\text{ if } i\leq r,\\
i-r, &\text{ if }i>r.\end{cases}$$
Following the construction we find the following realization of $ \hat{(\varphi(M^\ast))}$:
$$ \begin{pmatrix}
I_{n-r}&-C^\top & I_{n-r} & 0_r\\
(C & I_r & 0_{n-r}& I_r)^\epsilon
\end{pmatrix}.$$
By row permutation and multiplying columns with respective powers of $\epsilon$ we find
$$\begin{pmatrix}
C & I_r & 0_{n-r}&I_r\\
(I_{n-r}& -C^\top & I_{n-r} & 0_r)^\epsilon
\end{pmatrix}.$$
Hence, exchanging $A$ and $B$ yields a realization of $\hat{M}$. As a result, the face lattices correspond (with respect to $\varphi$ and exchanging $A$ and $B$) and, since $\varphi$ only acts on $E$, the claim follows.
\end{proof}

\subsubsection{Direct Sum}

As the Tutte polynomial of the direct sum of two matroids $M_1$ and $M_2$ factorizes, i.e. (see e.g., \cite{welshTutte})
$$T_{M_1\oplus M_2}(x,y)=T_{M_1}(x,y)\cdot T_{M_2}(x,y),$$
this is also true for our polynomial. For the proof, we make use of the following.\\

\begin{theorem}\label{directsum}
Let $M_1,M_2$ be two oriented matroids on disjoint groundsets $E_1,E_2$, respectively, and let $M_1\oplus M_2$ denote their direct sum. Then
$$\hat{(M_1 \oplus M_2)}=\hat{M_1}\oplus \hat{M_2},$$
where $\hat{(M_1 \oplus M_2)}$ is constructed using a basis $\beta=\beta_1\dot{\cup}\beta_2$ of $M_1\oplus M_2$ such that $\beta_1$ is a basis of $M_1$ and $\beta_2$ is a basis of $M_2$.
\end{theorem}
 
\begin{proof}
Let $B_1\dot{\cup} B_2$ be a basis of $\hat{(M_1\oplus M_2)}$, where $B_1$ is a basis of $(M_1\oplus M_2)_1$ and $B_2$ is a basis of $(M_1 \oplus M_2)_2$. Then $B_1=b_1\dot{\cup} b_2$ and $B_2=c_1\dot{\cup} c_2$, where $b_1$ is a basis of $(M_1)_1$, i.e., $M_1$ extended by parallels of the chosen basis $\beta_1$, $b_2$ is a basis of $(M_2)_1$, $c_1$ is a basis of $(M_1)_2$ and $c_2$ is a basis of $(M_2)_2$. Thus, 
$$ B_1\dot{\cup} B_2=(b_1\dot{\cup} b_2)\dot{\cup} (c_1\dot{\cup} c_2)=(b_1\dot{\cup} c_1)\dot{\cup} (b_2\dot{\cup} c_2)$$
and $b_1\dot{\cup} c_1$ is a basis of $\hat{M_1}$ and $b_2\dot{\cup} c_2$ is a basis of $\hat{M_2}$.\\
If $B_1\dot{\cup} B_2$ is chosen lexicographically minimal, then also $b_1\dot{\cup}c_1$ and $b_2\dot{\cup}c_2$ are lexicographically minimal. By Proposition 7.6.4 and Proposition 7.6.1(c),
\begin{flalign*}
\chi_{\hat{(M_1\oplus M_2)}}(B_1\dot{\cup}B_2)&=\chi_{(M_1\oplus M_2)_1}(B_1)\chi_{(M_1\oplus M_2)_2}(B_2)\\
&=\chi_{(M_1\oplus M_2)_1}(b_1\dot{\cup}b_2)\chi_{(M_1\oplus M_2)_2}(c_1\dot{\cup}c_2)\\
&=\chi_{(M_1)_1}(b_1)\chi_{(M_2)_1}(b_2)\chi_{(M_1)_2}(c_1)\chi_{(M_2)_2}(c_2)\\
&=\chi_{\hat{M_1}}(b_1\dot{\cup}c_1)\chi_{\hat{M_2}}(b_2\dot{\cup}c_2).
\end{flalign*}
\end{proof} 

\begin{lemma}\label{lemma}
\hfill
\begin{itemize}
\item[(i)] $\mathcal{L}(\hat{M_1}\oplus \hat{M_2})=\lbrace X_1\dot{\cup}X_2: X_1\in \mathcal{L}(\hat{M_1}),X_2\in \mathcal{L}(\hat{M_2})\rbrace$
\item[(ii)] $\rk_{\mathcal{L}(\hat{M_1}\oplus \hat{M_2})}(X)=\rk_{\mathcal{L}(\hat{M_1})}(X_1)+\rk_{\mathcal{L}(\hat{M_2})}(X_2)$, where $X_1 \in \mathcal{L}(\hat{M_1})$ and $X_2\in \mathcal{L}(\hat{M_2})$.
\item[(iii)] $\mu_{\mathcal{L}(\hat{M_1}\oplus \hat{M_2})}(\mathbf{0},X)=\mu_{\mathcal{L}(\hat{M_1})}(\mathbf{0},X_1)\cdot \mu_{\mathcal{L}(\hat{M_2})}(\mathbf{0},X_2)$, where $X_1 \in \mathcal{L}(\hat{M_1})$ and $X_2\in \mathcal{L}(\hat{M_2})$.

\end{itemize}
\end{lemma}
\begin{proof}
\begin{itemize}
\item[(i)] See Proposition 7.6.1(e) in \cite{BLSW}.
\item[(ii)] Follows directly from (i).
\item[(iii)] Lemma \ref{coneOM} and (ii) yield
\begin{flalign*}
\mu_{\mathcal{L}(\hat{M_1}\oplus \hat{M_2})}(\mathbf{0},X)&=\mu_{\mathcal{L}(\hat{M_1}\oplus \hat{M_2})}(\mathbf{0},X_1\dot{\cup}X_2)\\
&=(-1)^{\rk_{\mathcal{L}(\hat{M_1}\oplus \hat{M_2})}(X_1\dot{\cup}X_2)}\\
&=(-1)^{\rk_{\mathcal{L}(\hat{M_1})}(X_1)+\rk_{\mathcal{L}(\hat{M_2})}(X_2)}\\
&=(-1)^{\rk_{\mathcal{L}(\hat{M_1})}(X_1)}\cdot (-1)^{\rk_{\mathcal{L}(\hat{M_2})}(X_2)}\\
&=\mu_{\mathcal{L}(\hat{M_1})}(\mathbf{0},X_1)\cdot \mu_{\mathcal{L}(\hat{M_2})}(\mathbf{0},X_2),
\end{flalign*}
where $X_1 \in \mathcal{L}(\hat{M_1})$ and $X_2\in \mathcal{L}(\hat{M_2})$.
\end{itemize}
\end{proof}

\begin{corollary}
Let $M_1,M_2$ be two oriented matroids on disjoint groundsets $E_1,E_2$, respectively, and let $M_1\oplus M_2$ denote their direct sum. Let $\hat{(M_1 \oplus M_2)}$ be constructed using a basis $\beta=\beta_1\dot{\cup}\beta_2$ of $M_1\oplus M_2$ such that $\beta_1$ is a basis of $M_1$ and $\beta_2$ is a basis of $M_2$. Then
$$\Omega_{NL}^{M_1\oplus M_2, \beta}(x,y,z)=\Omega_{NL}^{M_1,\beta_1}(x,y,z)\cdot \Omega_{NL}^{M_2,\beta_2}(x,y,z).$$
\end{corollary}

\begin{proof}
By Theorem \ref{directsum} and Lemma \ref{lemma} we find
\begin{itemize}
\item[(i)] $\mathcal{L}_+(\hat{M_1\oplus M_2})=\lbrace X_1\dot{\cup} X_2: X_1\in \mathcal{L}_+(\hat{M_1}), X_2\in \mathcal{L}_+(\hat{M_2})\rbrace$
\item[(ii)] $\rk_{\mathcal{L}(\hat{M_1\oplus M_2})}(X)=\rk_{\mathcal{L}(\hat{M_1})}(X_1)+\rk_{\mathcal{L}(\hat{M_2})}(X_2)$, where $X_1 \in \mathcal{L}(\hat{M_1})$ and $X_2\in \mathcal{L}(\hat{M_2})$.
\item[(iii)] $\mu_{\mathcal{L}(\hat{M_1\oplus M_2})}(\mathbf{0},X)=\mu_{\mathcal{L}(\hat{M_1})}(\mathbf{0},X_1)\cdot \mu_{\mathcal{L}(\hat{M_2})}(\mathbf{0},X_2)$, where $X_1 \in \mathcal{L}(\hat{M_1})$ and $X_2\in \mathcal{L}(\hat{M_2})$.
\end{itemize}
Furthermore,
\begin{itemize}
\item[(iv)] $E\setminus (\underline{X}\cap E)= (E_1\dot{\cup}E_2)\setminus (\underline{X}\cap (E_1\dot{\cup}E_2))=(E_1\setminus( \underline{X}\cap E_1))\dot{\cup}( E_2\setminus (\underline{X}\cap E_2))$ 
\item[(v)] $\underline{X}\cap A= \underline{X}\cap (A_1\dot{\cup}A_2)=(\underline{X} \cap A_1)\dot{\cup}(\underline{X}\cap A_2)$
\item[(vi)]$\underline{X}\cap B= \underline{X}\cap (B_1\dot{\cup}B_2)=(\underline{X} \cap B_1)\dot{\cup}(\underline{X}\cap B_2)$
\end{itemize}
hold by construction, where in (v) and (vi) $A_1\subseteq A$ equals the set of elements parallel to those in $\beta_1$ and $A_2\subseteq A$ equals the set of elements parallel to those in $\beta_2$. Analogously, $B_1\subseteq B$ equals the set of elements parallel to those in $E_1\setminus \beta_1$ and $B_2\subseteq B$ equals the set of elements parallel to those in $E_2\setminus \beta_2$.\\
Alltogether,
\begin{flalign*}
\Omega_{NL}^{M_1\oplus M_2,\beta}&(x,y,z)=\sum_{X\in \hat{\mathcal{L}}_+}\hat{\mu}(\mathbf{0},X)x^{\rk_{\hat{\mathcal{L}}}(X)+|E\setminus (\underline{X}\cap E)|}y^{|\underline{X}\cap A|}z^{|\underline{X}\cap B|}\\
&=\sum_{X_1\in {\mathcal{L}}_+(\hat{M_1})}{\mu}_{\mathcal{L}(\hat{M_1})}(\mathbf{0},X_1)x^{\rk_{\hat{\mathcal{L}}(M_1)}(X_1)+|E_1\setminus (\underline{X_1}\cap E_1)|}y^{|\underline{X_1}\cap A_1|}z^{|\underline{X_1}\cap B_1|}\\
&\cdot \sum_{X_2\in {\mathcal{L}}_+(\hat{M_2})}{\mu}_{\mathcal{L}(\hat{M_2})}(\mathbf{0},X_2)x^{\rk_{\hat{\mathcal{L}}(M_2)}(X_2)+|E_2\setminus (\underline{X_2}\cap E_2)|}y^{|\underline{X_2}\cap A_2|}z^{|\underline{X_2}\cap B_2|} \\
&=\Omega_{NL}^{M_1,\beta_1}(x,y,z)\cdot \Omega_{NL}^{M_2,\beta_2}(x,y,z).
\end{flalign*}
\end{proof}

\section{Discussion}

Besides Observation \ref{111} we currently do not have any combinatorial interpretations of concrete evaluations of the dichromate. Some might be found in evaluations of $\Omega_{NL}^{M,\mathcal{B}}$ at $(x, 1, 1),(-1, 1, 0)$ and $ (-1, 0, 1).$\\

Furthermore, we are not aware of any meaningful interpretation in the non-regular case. Nevertheless the polynomial exists in this case and since the union does not need to preserve regularity we have in any event already crossed this line.\\

\begin{observation}\label{obs:regular}
Let $M$ be a regular oriented matroid of rank $r$ represented by $\begin{pmatrix}
I_r& C\end{pmatrix}$ such that $C$ is not identical to the zero matrix. Then $\hat{M}$ is not regular and in particular not graphic.
\end{observation}

\begin{proof}
Since $C$ is not the zero matrix there exists an entry $c_{ij}\neq 0$. Due to the construction of $\hat{M}$ we find a non-zero entry $-c_{ji}$ in $\hat{M}$ as well as a $1$ at position $(i,i)$ and at $(j,j)$. Together with the corresponding elements of the basis $A \cup B$ these entries build the four point line
$$ \begin{pmatrix}
1&c_{ij}&1&0\\-c_{ij}&1&0&1
\end{pmatrix}$$
as a minor. Thus, $\hat{M}$ can neither be regular nor graphic.
\end{proof}

Another beauty mistake is that the dichromate highly depends on the chosen basis.\\

\begin{example}
We examine the dichromate of the following orientation of the $K_4$ with respect to the chosen basis $\mathcal{B}$.
\begin{figure}[h!]
\begin{center}
\includegraphics[scale=0.6]{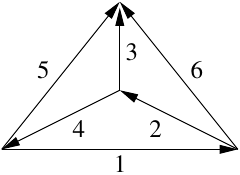}
\end{center}
\vspace{-0.5cm}
\end{figure}

\begin{itemize}
\item $\mathcal{B}=\lbrace 3,5,6\rbrace$ is a directed star
$$ M=\begin{pmatrix}
 1 & 0 & 0 & 0 & -1& 1\\
 0 & 1 & 0 & 1 & 0 & -1\\
 0 & 0 & 1 & -1& 1 &0
\end{pmatrix}$$ 
$$\Rightarrow \Omega_{NL}^{K_4,\mathcal{B}}(x,y,z)=x^6- x^4y^3 - x^4z^3 -2x^6y^3z^3 + 3x^5y^3z^2 + 3x^5y^2z^3 - 3x^4y^2z^2 $$

\item $\mathcal{B}=\lbrace 1,4,5\rbrace$ is an undirected star
$$M=\begin{pmatrix}
1 & 0 & 0 & -1  & 0 & -1\\
 0 & 1 & 0 & -1 & 1 &0\\
 0 & 0 & 1 & 0  & 1 &1
\end{pmatrix}$$ 
$$\Rightarrow \Omega_{NL}^{K_4,\mathcal{B}}(x,y,z)=x^6 - x^4y - x^4z + x^2yz$$
\item $\mathcal{B}=\lbrace 1,2,3\rbrace$ is a directed path
$$M=\begin{pmatrix}
 1 & 0 & 0 & -1 & 1 & 0\\
 0 & 1 & 0 & -1&  1 & 1\\
 0 & 0 & 1 & 0  & 1&  1
\end{pmatrix}$$ 
$$\Rightarrow \Omega_{NL}^{K_4,\mathcal{B}}(x,y,z)=x^6 - x^4y - x^4z + x^2yz$$
\item $\mathcal{B}=\lbrace 1,3,4\rbrace$ is an undirected path
$$M=\begin{pmatrix}
 1 & 0 & 0 & -1& 0 & -1\\
 0 & 1 & 0 & 0 & 1 &1\\
 0 & 0 & 1 &-1 &-1 & -1
\end{pmatrix}$$ 
$$\Rightarrow \Omega_{NL}^{K_4,\mathcal{B}}(x,y,z)=x^6 - x^4y - x^4z + x^2yz$$
\item $\mathcal{B}=\lbrace 1,3,6\rbrace$ is another undirected path
$$M=\begin{pmatrix}
 1 & 0 & 0 & 0& -1 & 1\\
 0 & 1 & 0 & -1 & 1 &0\\
 0 & 0 & 1 &1 &-1 & 1
\end{pmatrix}$$ 
$$\Rightarrow \Omega_{NL}^{K_4,\mathcal{B}}(x,y,z)=x^6 + x^4y^2z^2 - x^4y^2 - x^4z^2$$
\end{itemize}
That the polynomials coincide for the bases $\lbrace 1,2,3\rbrace$ and $\lbrace 1,3,4\rbrace$ is due to Theorem \ref{triduality}. Presently we are not aware of an explanation why the polynomial for $\lbrace 1,4,5\rbrace$ also coincides with those.
\end{example}

The dependence on the basis might be resolved (or at least better understood) by looking at (directed) internal and external activities (see \cite{tuttecon}).\\

Since the contraction of arcs might generate new directed cycles and loops it is clear that our polynomials do not satisfy the (classical) deletion-contraction formula. In \cite{moreno:recursive} a recursive formula is given, using a different definition of contraction. Presumably the most agreed concept of digraph minors in the context of acyclic colorings are \emph{butterfly minors} (see \cite{johnson}). Unfortunately digraphs that are not butterfly contractible can be arbitrarily complicated. \\

\bibliography{NL-flow}{}

\begin{thebibliography}{10}

\bibitem{althowie}
{\sc B.~Altenbokum, W.~Hochst\"attler, and J.~Wiehe}, {\em The {NL}-flow
  polynomial}, Discrete Applied Mathematics, 296 (2021), pp.~193--202.
\newblock 16th Cologne-Twente Workshop on Graphs and Combinatorial Optimization
  (CTW 2018).

\bibitem{BLSW}
{\sc A.~Bj\"orner, M.~Las~Vergnas, B.~Sturmfels, N.~White, and G.~M. Ziegler},
  {\em Oriented Matroids}, vol.~46 of Encyclopedia of Mathematics and its
  Applications, Cambridge University Press, 2~ed., 1999.

\bibitem{bondy}
{\sc J.~A. Bondy and U.~S.~R. Murty}, {\em Graph Theory}, Springer London,
  2008.

\bibitem{dualpairs}
{\sc T.~H. Brylawski and G.~M. Ziegler}, {\em Topological representation of
  dual pairs of oriented matroids}, Discrete \& Computational Geometry, 10
  (1993), pp.~237--240.

\bibitem{moreno:recursive}
{\sc D.~Gonz\'{a}lez-Moreno, R.~Hern\'{a}ndez-Ortiz, B.~Llano, and M.~Olsen},
  {\em The dichromatic polynomial of a digraph}, Graphs and Combinatorics, 38
  (2022).

\bibitem{hochstarticle}
{\sc W.~Hochst{\"a}ttler}, {\em A flow theory for the dichromatic number},
  European Journal of Combinatorics, 66 (2017), pp.~160--167.

\bibitem{howie}
{\sc W.~Hochst{\"a}ttler and J.~Wiehe}, {\em The chromatic polynomial of a
  digraph}, in Graphs and Combinatorial Optimization: from Theory to
  Applications, CTW2020 Proceedings, 2021, pp.~1--14.

\bibitem{johnson}
{\sc T.~Johnson, N.~Robertson, P.~D. Seymour, and R.~Thomas}, {\em Directed
  tree-width}, Journal of Combinatorial Theory. Series B,  (2001),
  pp.~138--154.

\bibitem{neumannlara}
{\sc V.~Neumann-Lara}, {\em The dichromatic number of a digraph}, Journal of
  Combinatorial Theory. Series B, 33 (1982), pp.~265--270.

\bibitem{tuttecon}
{\sc W.~T. Tutte}, {\em A contribution to the theory of chromatic polynomials},
  Canadian Journal of Mathematics, 6 (1954), pp.~80--91.

\bibitem{welshTutte}
{\sc D.~J.~A. Welsh}, {\em The {T}utte polynomial}, Random Structures \&
  Algorithms, 15 (1999), pp.~210--228.

\end{thebibliography}
\bibliographystyle{siam}

\end{document}